\documentclass{article}

\usepackage{amssymb}
\usepackage{latexsym}
\usepackage{epsfig}

\newtheorem{iresult}{Theorem}

\newtheorem{theorem}{Theorem}[section]
\newtheorem{lemma}[theorem]{Lemma}
\newtheorem{corollary}[theorem]{Corollary}
\newtheorem{proposition}[theorem]{Proposition}

\newenvironment{definition}{\refstepcounter{theorem}\smallskip\par\noindent\textbf{Definition\hspace{0.7ex}
\thetheorem}\,}{\smallskip\par\noindent}
\newenvironment{remark}{\refstepcounter{theorem}\smallskip\par\noindent\textbf{Remark\hspace{0.7ex}
\thetheorem}\,}{\smallskip\par\noindent}
\newenvironment{example}{\refstepcounter{theorem}\smallskip\par\noindent\textbf{Example\hspace{0.7ex}
\thetheorem}\,}{\smallskip\par\noindent}

\newcommand{\proof}{\par\noindent\textit{Proof.}\,\,}
\newcommand{\eproof}{\hspace*{\fill}$\Box$\smallskip\par\noindent}

\newcommand{\bL}{\mathbf{L}}

\newcommand{\cI}{\mathcal{I}}
\newcommand{\cJ}{\mathcal{J}}

\newcommand{\cM}{\mathcal{M}}
\newcommand{\cO}{\mathcal{O}}
\newcommand{\cP}{\mathcal{P}}

\newcommand{\cX}{\mathcal{X}}

\newcommand{\cC}{\mathcal{C}}

\newcommand{\Proj}{\mathrm{Proj}\,}

\newcommand{\Reg}{\mathrm{Reg}\,}

\newcommand{\Frac}{\mathrm{Frac}\,}

\newcommand{\Ker}{\mathrm{Ker}\,}
\newcommand{\Rad}{\mathrm{Rad}\,}

\newcommand{\pdim}{\mathrm{pdim}\,}
\newcommand{\fdim}{\mathrm{fdim}\,}

\newcommand{\wdim}{\mathrm{wdim}\,}
\newcommand{\height}{\mathrm{height}\,}
\newcommand{\gldim}{\mathrm{gldim}\,}

\newcommand{\grade}{\mathrm{grade}\,}
\newcommand{\Grade}{\mathrm{Grade}\,}
\newcommand{\trdeg}{\mathrm{trdeg}\,}
\newcommand{\ch}{\mathrm{char}\,}

\newcommand{\Spec}{\mathrm{Spec}\,}

\newcommand{\fgSpec}{\mathrm{fgSpec}\,}
\newcommand{\MaxSpec}{\mathrm{MaxSpec}\,}

\newcommand{\sep}{\;|\;}

\begin{document}
\begin{titlepage}
\title{Regular local algebras over a Pr\"ufer domain:\\
weak dimension and regular sequences}
\author{Hagen Knaf\\
Fraunhofer Institut Techno- und Wirtschaftsmathematik\\
Gottlieb-Daimler-Strasse 49,  67663 Kaiserslautern, Germany\\
email: knaf@itwm.fhg.de}
\end{titlepage}
\maketitle

 %Inserted by TeXtelmExtel

\begin{abstract}
A not necessarily noetherian local ring $O$ is called regular if every finitely generated ideal $I\lhd O$ possesses finite projective dimension. In the article localizations $O=A_q$, $q\in\Spec A$, of a finitely presented, flat algebra $A$ over a Pr\"ufer domain $R$ are investigated with respect to regularity: this property of $O$ is shown to be equivalent to the finiteness of the weak homological dimension $\wdim O$. A formula to compute $\wdim O$ is provided. Furthermore regular sequences within the maximal ideal $M\lhd O$ are studied: it is shown that regularity of $O$ implies the existence of a maximal regular sequence of length $\wdim O$. If $\height (q\cap R)\neq\infty$, then this sequence can be choosen such that the radical  of the ideal generated by the members of the sequence equals $M$. As a consequence it is proved that if $O$ is regular, then the (noetherian) factor ring $O/(q\cap R)O$ is Cohen-Macaulay. If $(q\cap R)R_{q\cap R}$ is not finitely generated, then $O/(q\cap R)O$ itself is regular.
\end{abstract}

 %Inserted by TeXtelmExtel

\newpage
\section*{Introduction}
In the present article finitely generated, flat algebras $A$ over a Pr\"ufer domain $R$ and  their localizations $A_q$ at prime ideals $q\lhd A$ are investigated. The finiteness properties and the ideal theory of these algebras have been the subject of various investigations, two of the main results being:
\begin{itemize}
\item coherence \cite{Gla1}, Ch.~7: every finitely generated ideal $I\lhd A$ is finitely presented,
\item catenarity \cite{BDF}: the length $\ell$ of a non-refinable chain of primes $q_0\subset q_1\subset\ldots\subset q_\ell$ in $A$ depends on $q_0$ and $q_\ell$ only, provided that the Krull dimension $\dim R_p$ is finite for all $p\in\Spec R$.
\end{itemize}
Motivated by questions concerning the geometry of schemes over Pr\"ufer domains in the sequel we study localizations $A_q$ possessing a distinguished property that generalizes the noetherian notion of regularity. A noetherian local ring $O$ is called regular if its global homological dimension $\gldim O$ is finite. This property plays an important role in the theory of noetherian rings and links algebra with geometry. Accordingly attempts were made to generalize it to non-noetherian rings. They led to the following definition first given by J.~Bertin \cite{Ber}: the local ring $O$ is called regular if for every {\em finitely generated} ideal $I\lhd O$ the projective dimension $\pdim_O O/I$ of the $O$-module $O/I$ is finite. Auslander's Lemma $\gldim O=\sup (\pdim_O (O/I)\sep I\lhd O)$ shows that Bertin's definition indeed generalizes the noetherian notion of regularity. However a local ring regular in the sense of Bertin needs not have finite global dimension.

 %Inserted by TeXtelmExtel

In the sequel the words {\em regular} and {\em regularity} always refer to Bertin's definition.

 %Inserted by TeXtelmExtel

The weak homological dimension $\wdim O$ of a local ring $O$ in many cases has shown to be an appropriate tool to check for regularity. Analogously to the global dimension it is defined as the supremum of the flat dimensions $\fdim M$ of all $O$-modules $M$, where $\fdim M$ denotes the length of the shortest flat resolution of $M$.  It is well-known that a coherent local ring $O$ of finite weak dimension $\wdim O$ is regular and that the finiteness of the global dimension $\gldim O$ implies that of $\wdim O$. However there exist coherent, regular, local rings possessing an infinite weak dimension, so that in contrast to the noetherian case regularity of a coherent local ring cannot in general be checked using the weak dimension. 

 %Inserted by TeXtelmExtel

Several classes of coherent local rings for which regularity is {\em equivalent} to the finiteness of the weak dimension have been identified in the past -- an overview can be found in \cite{Gla1}. Among them are for example the local rings with finitely generated maximal ideal \cite{TZT} and localizations of certain group algebras. One of the main results of the present article describes a new class of local rings of that type: 
\newpage
\begin{iresult}
\label{main 1}
The localization  $A_q$, $q\in\Spec A$, of a finitely generated, flat algebra $A$ over the Pr\"ufer domain $R$ is regular if and only if its weak dimension is finite. For a regular localization $A_q$ the  weak dimension satisfies:
\[
\wdim A_q=\left\{
\begin{array}{ll}
\dim (A_q\otimes_R kp)&\mbox{ if $p=0$}\\
\dim (A_q\otimes_R kp)+1&\mbox{ if $p\neq 0$},
\end{array}
\right.
\]
where $p:=q\cap R$ and $kp$ denotes the field of fractions of $R/p$.
\end{iresult}
As the second theme of this article we discuss a method to construct regular sequences within a localization $A_q$. It is copied from the noetherian case and is shown to work for an arbitrary coherent local ring of finite weak dimension:
\begin{iresult}
\label{main 2}
In a coherent local ring $O$ of finite weak dimension every sequence $(t_1,\ldots ,t_\ell)$ composed of elements of the maximal ideal $M\lhd O$ such that $t_1+M^2,\ldots ,t_\ell +M^2\in M/M^2$ are $O/M$-linearly independent is a regular sequence. Its members generate a prime ideal $q\lhd O$ with the property 
\[
\wdim O/q=\wdim O-\ell ;
\]
in particular $O/q$ is a coherent regular ring. 
\end{iresult}
It is well-known that under the assumptions made in Theorem \ref{main 2} the length $\ell$ of a regular sequence within the maximal ideal satisfies $\ell\leq\wdim O$. The theorem thus yields the inequality
\begin{equation}
\label{dimension bound}
\dim M/M^2\leq\wdim O.
\end{equation}
In contrast to the noetherian case however lifting a basis of $M/M^2$ does not in general yield a {\em maximal} regular sequence. Moreover the length of the longest regular sequence in $M$  needs not be equal to $\wdim O$. Against this background the behavior of the lifting procedure in the case of algebras over Pr\"ufer domains is remarkable:
\begin{iresult}
\label{main 3}
The maximal ideal $M:=qA_q$ of a regular localization $O:=A_q$ of a finitely generated, flat algebra $A$ over a Pr\"ufer domain $R$ contains a maximal regular sequence $(t_1,\ldots ,t_d)$ of length $d=\wdim O$.

 %Inserted by TeXtelmExtel

Every sequence $(t_1,\ldots ,t_\ell)$ of elements of $M$ such that  $(t_1+M^2,\ldots ,t_\ell+M^2)$ forms an $O/M$-basis of $M/M^2$ is regular and has the following properties:
\newpage
\begin{enumerate}
\item If $M$ is finitely generated, then $(t_1,\ldots ,t_\ell )$ is a maximal regular sequence of length $\ell=\wdim O$ and its members generate $M$.
\item If $M$ is not finitely generated, then for every non-zero $t\in p:=M\cap R$ the sequence $(t_1,\ldots ,t_\ell,t)$ is maximal regular; its length satisfies $\ell +1=\wdim O$. Moreover if $pR_p$ is the radical of the principal ideal $tR_p$, then $M$ equals the radical of the ideal $\sum\limits_{i=1}^\ell Ot_i +Ot$.
\end{enumerate}
\end{iresult}
\noindent
Note that as a consequence of Theorem \ref{main 1} and inequality (\ref{dimension bound}) the dimension of $M/M^2$ in Theorem \ref{main 3} is finite. Furthermore notice that the maximal ideal $M$ is finitely generated if the maximal ideal $pR_p$, $p=M\cap R$, is finitely generated. Finally it should be mentioned that the local ring $R_p$ is a valuation domain. Thus if $\dim R_p$ is finite, then $pR_p$ is the radical of a principal ideal.
\medskip
 
 %Inserted by TeXtelmExtel
 
The results presented so far can be used to get insight into the geometry of integral separated $R$-schemes $\cX$ of finite type over a Pr\"ufer domain $R$ -- in the sequel such schemes are called $R$-varieties. They are known to be $R$-flat and to possess coherent structure sheaf $\cO_\cX$. Furthermore all fibres $\cX\times_Rkp$, $p\in\Spec R$, are equidimensional schemes of a common dimension $n\in\mathbb{N}$; in particular the generic fibre $X:=\cX\times_RK$, $K:=\Frac R$, is an algebraic variety of dimension $n$. The case of a noetherian base ring $R$, that is $R$-varieties over a Dedekind domain, are studied extensively in arithmetic algebraic geometry. $R$-varieties over a non-noetherian valuation domain appear in the valuation theory of algebraic function fields: the function field $F$ of a {\em normal} $R$-variety $\cX$ carries a finite set $V$ of non-archimedian valuations induced by $\cX$. Namely the valuation ring of each of the valuations $v\in V$ is one of the local rings $\cO_{\cX ,\eta}$, where $\eta\in\cX$ runs through the generic points of the closed fibre $\overline{X}:=\cX\times_R kN$, $N\lhd R$ the maximal ideal of $R$. By definition the restriction $v|_K$ is a valuation $v_K$ of $K$ independent of $v$ and possessing $R$ as its valuation domain. Moreover the extension of residue fields of $v|v_K$ equals $k\eta |kN$ and thus has the same transcendence degree as $F|K$. Valuations with the latter property are called constant reductions of $F|K$. They were first studied using a purely valuation-theoretic approach in the case of transcendence degree one by M. Deuring \cite{Deu} and E. Lamprecht \cite{Lam}. For higher transcendence degree P. Roquette \cite{Roq} introduced a geometric approach by (essentially) studying projective $R$-varieties. His work supplements and generalizes results obtained by G. Shimura \cite{Shi} in the case of a discrete valuation ring $R$.

 %Inserted by TeXtelmExtel

In the articles just mentioned a focus lies on clarifying the relation between the divisor theory of the generic fibre $X$ of an $R$-variety $\cX$ and that of its closed fibre $\overline{X}$. In this task Weil divisors on $\cX$ itself enter in. On the non-noetherian scheme $\cX$ they must be defined utilizing sheaves instead of taking the usual approach through subschemes of codimension one: a Weil divisor on $\cX$ is a coherent fractional  $\cO_{\cX}$-module $\cJ$ with the property $\cJ=(\cO_\cX :(\cO_\cX :\cJ ))=:\widehat{\cJ}$ called reflexivity. The set $D(\cX )$ of all Weil divisors on $\cX$ becomes an abelian semigroup through the multiplication $\cI\circ\cJ:=\widehat{\cI\cJ}$. If $\cX$ is normal, then $D(\cX )$ is a group \cite{Kn1}, that in the noetherian case is isomorphic to the group of ordinary Weil divisors. Given a Weil prime divisor $P\subset X$ of the generic fibre of $\cX$ a divisor $\overline{P}$ on $\overline{X}$ can be defined essentially by intersecting the Zariski closure $\cP$ of $P$ on $\cX$ with the subscheme $\overline{X}\subset\cX$. The ideal sheaf $\cJ$ associated to the subscheme $\cP$ is an element of $D(\cX )$ -- a non-trivial result in the present non-noetherian context. However to define a divisor on the possibly singular scheme $\overline{X}$ it is desirable that $\cJ$ be invertible. At that point regularity becomes important. 
 
 %Inserted by TeXtelmExtel
 
Every stalk $\cJ_P$, $P\in\cX$, of a sheaf $\cJ\in D(\cX )$ is a fractional $\cO_{\cX ,P}$-ideal with the property 
$(\cO_{\cX ,P}:(\cO_{\cX ,P}:\cJ_P ))=\cJ_P$. Such fractional ideals are principal if in the local ring $\cO_{\cX ,P}$ finitely many elements always possess a greatest common divisor. One of the main results of the theory of finite free resolutions states that a coherent, regular, local ring is of that type -- see \cite{Nor}. Thus as in the noetherian case the sheaves $\cJ\in D(\cX )$ are locally principal on the regular locus $\Reg\cX :=\{P\in\cX\sep\cO_{\cX ,P}\mbox{ is regular.}\}$ and on a regular $R$-variety -- meaning that $\Reg\cX =\cX$ -- Weil and Cartier divisors coincide.

 %Inserted by TeXtelmExtel

A description of the regular locus can be achieved using Theorem \ref{main 1}:
\smallskip

 %Inserted by TeXtelmExtel

\noindent
\textit{The regular locus of an $R$-variety $\cX$ of relative dimension $n\in\mathbb{N}$ can be characterized as $\Reg\cX =\{P\in\cX | \wdim\cO_{\cX,P}\leq n+1\}$ and satisfies
\[
\bigcup\limits_{p\in\Spec R}\Reg (\cX\times_Rkp)\subseteq\Reg\cX.
\]}
A regular algebraic variety $X$ over the field $K$ is said to have good reduction at the valuation ring $R\subset K$ if there exists an $R$-variety $\cX$ such that $X=\cX\times_RK$ and the closed fibre $\cX\times_RkN$ is regular. If $R$ is noetherian, then the scheme $\cX$ itself is regular. Vice versa it is known that the closed fibre  $\cX\times_RkN$ of a regular variety over a discrete valuation domain $R$ is Cohen-Macaulay. As the main geometric result of the present article we prove that these facts remain valid over a general Pr\"ufer domain $R$ -- with an interesting non-noetherian twist though:
\begin{iresult}
\label{main4}
The fibres $\cX\times_Rkp$ of a regular $R$-variety $\cX$ over the Pr\"ufer domain $R$ are Cohen-Macaulay. If the prime ideal $pR_p$ is not finitely generated, then the fibre $\cX\times_Rkp$ is a regular $kp$-scheme.

 %Inserted by TeXtelmExtel

Vice versa: an $R$-variety $\cX$ with closed structure morphism $\cX\rightarrow\Spec R$ and such that all fibres $\cX\times_Rkp$ over maximal ideals $p\in\Spec R$ are regular is itself regular.
\end{iresult}
As an end of this introduction some facts about normal {\em curves} over a non-noetherian valuation domain $R$ should  be mentioned: first of all the relationship between valued function fields of transcendence degree one and proper normal $R$-curves was clarified by B. Green, M. Matignon and F. Pop \cite{GMP3}, \cite{G} resulting in a category equivalence that generalizes the well-known equivalence between projective, normal, algebraic curves and function fields of transcendence degree one. This and related results found applications in the context of Rumely's Local-Global Principle \cite{GPR} and the so-called Skolem problems \cite{GMP4}. Earlier E.~Kani in his non-standard approach to Diophantine Geometry \cite{Kan1} established an intersection product on a valued function field making use of the regularity of $R$-curves of the form $C\times_k R$, where $C$ is a smooth curve over the field $k\subset R$.  
\smallskip

 %Inserted by TeXtelmExtel

\noindent
\textbf{Guide for the reader}

 %Inserted by TeXtelmExtel

Section 1 deals with regularity and the finiteness of the weak dimension; the proof of Theorem \ref{main 1} is provided. Moreover a brief introduction to non-noetherian regularity as well as to the basic properties of finitely generated algebras over Pr\"ufer domains is given.

 %Inserted by TeXtelmExtel

In Section 2 the construction of regular sequences through lifting linearly independent elements from $M/M^2$ to $M$, where $M$ is the maximal ideal of a coherent local ring of finite weak dimension, is investigated. The results are used to prove Theorem 2 and 3 as well as the first part of Theorem 4.

 %Inserted by TeXtelmExtel

Section 3 contains a geometric point of view of the results obtained in the Sections 1 and 2. A curve over a two-dimensional valuation domain is discussed as an example.

 %Inserted by TeXtelmExtel

In the appendix an upper bound for the global dimension $\gldim A_q$ of a regular localization of a finitely generated, flat $R$-algebra $A$ over a Pr\"ufer domain $R$ of small cardinality is given.
\smallskip

 %Inserted by TeXtelmExtel

\noindent
\textbf{Conventions:} Throughout the article a \textbf{local ring} is a commutative not necessarily noetherian ring $O$ with unity possessing a unique maximal ideal $M$. A local ring $(O,M)$ is said to be \textbf{essentially finitely generated} respectively \textbf{essentially finitely presented} over the ring $R$ if $O=A_q$, $q\in\Spec A$, for some finitely generated resp. finitely presented $R$-algebra $A$. Extensions $A\subseteq B$ of rings are denoted by $B|A$. If $A$ and $B$ are integral domains $\trdeg (B|A)$ denotes the transcendence degree of the field extension $\Frac B|\Frac A$. For an ideal $I\lhd A$ of a commutative ring $A$ the radical of $I$ is denoted by $\Rad I$.
\smallskip

 %Inserted by TeXtelmExtel

\noindent
\textbf{Acknowledgement}: the results presented in this article were partially worked out during two visits of the Department of Mathematics and Statistics, University of Saskatchewan at Saskatoon. I owe thanks to Salma and Franz Viktor Kuhlmann for the inspiring atmosphere they created during these visits and for financial support.
% ===================
\section{Weak dimension}
\label{s:wdim}
% ===================
Let $A$ be a finitely generated, flat algebra over a Pr\"ufer domain $R$. The main objective in this section is to determine the weak homological dimension of the localizations $O:=A_q$, $q\in\Spec A$. We are interested in the case of a non-noetherian domain $R$, so that in general the local ring $O$ is non-noetherian too. We shall see that the weak dimension of  $O$ is finite if and only if $O$ is regular in the sense of Bertin \cite{Ber}. In the regular case we express the weak dimension in terms of the dimension of the fibre ring $O\otimes_R kp$, $p:=q\cap R$. Bertin's generalization of the noetherian notion of regularity is central within the present article. We therefore give a brief summary of relevant facts in Subsection \ref{ss:regularity}. Readers acquainted with this topic may directly jump to Subsection \ref{ss:main} that contains the principal result of the present section, the determination of the weak dimension $\wdim O$. Auxiliary results of independent interest on upper bounds for the weak dimension are treated in Subsection \ref{ss:bounds}, while in the last Subsection \ref{ss:pruefer} basic properties of algebras over a Pr\"ufer domain are reviewed.
% ===========================
\subsection{Non-noetherian regularity}
\label{ss:regularity}
% ===========================
The local ring $O$ is called \textbf{regular} if the projective dimension $\pdim I$ of every {\em finitely generated} ideal $I\lhd O$ is finite. This property in general does {\em not} imply the finiteness of the \textbf{global dimension}
\[
\gldim O:=\sup (\pdim V\sep V\mbox{ an $O$-module})
\]
as in the noetherian case. Indeed every valuation domain is regular, however valuation domains can have infinite global dimension \cite{Os1}.

 %Inserted by TeXtelmExtel

If $O$ is a \textbf{coherent ring}, that is if every finitely generated ideal $I\lhd O$ is finitely presented, then flat resolutions 
can be used to check for regularity: the \textbf{flat dimension $\fdim V$} of an $O$-module $V$  is defined to be the length $n\in\mathbb{N}$ of the shortest resolution
\begin{equation}
\label{flat resolution}
0\rightarrow F_n\rightarrow ... \rightarrow F_0\rightarrow V\rightarrow 0
\end{equation}
of $V$ by flat $O$-modules $F_i$. If no such $n$ exists, then $\fdim V$ is set equal to $\infty$. The \textbf{weak dimension $\wdim O$} of $O$ is defined to be 
\[
\wdim O:=\sup (\fdim V\sep V\mbox{ an $O$-module}).
\]
Let $M$ be the maximal ideal of $O$ and consider the residue field $k:=O/M$ as an $O$-module. If $O$ is coherent, then we have (\cite{Gla1}, Thm.~3.1.3)
\begin{equation}
\label{wdim O=fdim k}
\wdim O=\fdim k.
\end{equation}
Moreover for a finitely generated ideal $I\lhd O$ the module $F_0$ in the resolution (\ref{flat resolution}) of $I$ can be choosen to be finitely generated. The coherence of $O$ then implies that all of the modules $F_i$ are finitely presented and thus projective. This shows: \textit{a coherent local ring $O$ of finite weak dimension is regular.}

 %Inserted by TeXtelmExtel

Unfortunately the class of coherent local rings of finite weak dimension does not contain all coherent, regular, local rings:
\begin{example}
The localization 
\[
O:=K[X_i\sep i\in\mathbb{N}]_P,\; P:=\sum\limits_{i\in\mathbb{N}}K[X_i\sep i\in\mathbb{N}]X_i ,
\]
of the polynomial ring in countably many variables over a field $K$ is a coherent, regular, local ring with the property $\wdim O=\infty$ -- see \cite{Gla1}, page 202.
\end{example}%
Note that in general the inequality $\wdim O\leq\gldim O$ holds. Strict inequality is widespread for non-noetherian rings $O$; in the noetherian case  of course equality is present.

 %Inserted by TeXtelmExtel

The formula (\ref{wdim O=fdim k}) frequently is not suitable to determine $\wdim O$. Instead as in the noetherian case regular sequences can be used as a device, although the theory of regular sequences in non-noetherian rings is far more complicated than in the noetherian case: let $I\lhd A$ be an ideal in the commutative ring $A$. Following Northcott \cite{Nor} the \textbf{(classical) grade of $I$} is defined to be the supremum of the lengths $\ell\in\mathbb{N}$ of the finite $A$-regular sequences $(a_1,\ldots ,a_\ell )$ contained in $I$; it is denoted by $\grade I$ and can attain the value $\infty$.

 %Inserted by TeXtelmExtel

The function $g(I,n):=\grade IA[X_1,\ldots ,X_n]$ assigning to each $n\in\mathbb{N}$ the grade of the lifting of the ideal $I\lhd A$ into the polynomial ring $A[X_1,\ldots ,X_n]$ is monotonically increasing. The limit
\[
\Grade I:=\lim_{n\to\infty}\grade IA[X_1,\ldots ,X_n]\in\mathbb{N}\cup\{\infty\}
\]
is called \textbf{polynomial} or \textbf{non-noetherian grade of $I$}. In the present context it can be used to determine weak dimensions (see \cite{Gla2}, Lemma 3): \textit{the weak dimension of a coherent, regular, local ring $O$ with maximal ideal $M$ satisfies
\begin{equation}
\label{wdim=Grade}
\wdim O=\Grade M.
\end{equation}
}
\noindent
In particular $\grade M\leq\wdim O$ holds in every coherent, regular, local ring, since $\grade I\leq\Grade I$ holds for every ideal $I\lhd O$.
% =============================
\subsection{A formula for the weak dimension}
\label{ss:main}
% =============================
Throughout the whole subsection let $R$ be a Pr\"ufer domain with field of fractions $K\neq R$. We consider localizations $O=A_q$ of a flat $R$-algebra $A$ at a prime $q\lhd A$. The local ring $S:=R_p$, $p:=q\cap R$, then is a valuation domain and the ring extension $O|S$ is faithfully flat. We are interested in calculating the weak dimension $\wdim O$ in the case where $A$ can be choosen to be finitely generated, that is the extension $O|S$ is essentially of finite type. The tools applied in this case yield an upper bound for $\wdim O$ of a regular ring $O$ belonging to a more general class of local $S$-algebras, namely:
\smallskip

 %Inserted by TeXtelmExtel

\noindent
($\bL$)
\textit{
For a valuation domain $S$ with maximal ideal $N$ let $\bL (S)$ be the class of local, faithfully flat, coherent $S$-algebras $O$ such that the maximal ideal $\overline{M}$ of the local ring $\overline{O}:=O/NO$ is finitely generated.
}
\smallskip

 %Inserted by TeXtelmExtel

\noindent
In the sequel we use the brief notation $(O,M)$ to denote a local ring $O$ having maximal ideal $M$. Consequently we sometimes write $\bL((S,N))$ for the class $\bL(S)$ defined in (L). We are always assuming that $S\neq\Frac S$ and call such a valuation domain \textbf{non-trivial}.

 %Inserted by TeXtelmExtel

In Subsection \ref{ss:pruefer} we shall prove that a finitely generated, flat $R$-algebra $A$ is finitely presented and coherent -- see Theorem \ref{coherence} and its corollary. Consequently every localization $A_q$, $q\in\Spec A$, is a member of the class $\bL (R_{q\cap R})$. 
The main result of the present section allows the calculation of the weak dimension of $A_q$: 
\begin{theorem}
\label{dimension equality}
Let $S$ be a non-trivial valuation domain. A essentially finitely presented $S$-algebra $O\in\bL (S)$ is regular if and only if the weak dimension $\wdim O$ is finite. If $O$ is regular, then $\wdim O=\dim\overline{O}+1$ holds.
\end{theorem}
\noindent
The proof of Theorem \ref{dimension equality} consists of three major steps:

 %Inserted by TeXtelmExtel

First an upper bound for the weak dimension in the regular case is provided. This result is valid for a larger class of local $S$-algebras than actually considered in Theorem \ref{dimension equality}. It involves the Zariski cotangent space of the local ring $(\overline{O},\overline{M})$ defined in (L): the $\overline{O}$-module $\overline{T}:=\overline{M}/\overline{M}^2$ is a vector space over the residue field $k=O/M=\overline{O}/\overline{M}$. Its dimension is finite by assumption and yields an upper bound for $\wdim O$:
\begin{theorem}
\label{upper bounds}
Let $S$ be a non-trivial valuation domain. The weak dimension of a regular local ring $O\in\bL (S)$ satisfies $\wdim O\leq\dim_k \overline{T}+1$, $k=\overline{O}/\overline{M}$. If the local ring $\overline{O}$ is noetherian, then the stronger inequality $\wdim O\leq\dim\overline{O}+1$ holds. 
In particular: a local ring $O\in\bL (S)$ is regular if and only if its weak dimension is finite. 
\end{theorem}
Theorem \ref{upper bounds} lays the basis for a inductive proof of Theorem \ref{dimension equality} once we are able to reduce the weak dimension of $O$ by a natural method. Such a method is provided by the following result well-known in the noetherian case: 
\begin{theorem}
\label{regular prime}
Let $(O,M)$ be a coherent, regular, local ring of finite weak dimension. Then for every element $t\in M\setminus M^2$ the local ring $O/tO$ is regular (and coherent) and the formula $\wdim (O/tO)=\wdim O-1$ holds. In particular $tO\in\Spec O$.
\end{theorem}
Finally we have to verify the assertion of Theorem \ref{dimension equality} for the starting point of the induction, that lies at $\wdim O=1$: for $O\in\bL (S)$ the extension $O|S$ is faithfully flat, thus by \cite{Gla1}, Thm.~3.1.1 $\wdim O\geq\wdim S$. On the other hand a local ring $S$ is a non-trivial valuation domain if and only if $\wdim S=1$  -- see \cite{Gla1}, Cor.~4.2.6.
\begin{lemma}
\label{wdim O=1}
Let $(S,N)$ be a non-trivial valuation domain. An essentially finitely presented $S$-algebra $(O,M)\in\bL (S)$ of weak dimension $\wdim O=1$ is a valuation domain with the property  $\Rad (NO)=M$.
\end{lemma}
\proof
We already mentioned that a local ring of weak dimension $1$ is a valuation domain. Choose some finitely presented, flat $S$-algebra $A$ such that $O=A_q$ for some $q\in\Spec A$ lying over $N$. In Subsection \ref{ss:pruefer}, Lemma \ref{domain} we shall see that since $O$ is a domain $A$ can be choosen to be a domain too. An application of \cite{Kn1}, Lemma 2.8 yields that the prime $M\cap A$ is minimal among the primes of $A$ lying over $N$. Since the primes contained in $M\cap A$ are totally ordered with respect to inclusion, the minimality of $M$ implies $\Rad (NO)=M$.
\eproof
We now directly turn  to the proof of Theorem \ref{dimension equality} postponing the proofs of Theorem \ref{upper bounds} and \ref{regular prime} to the Subsections \ref{ss:bounds} and  \ref{ss:regseq+wdim} respectively. The required properties of algebras over Pr\"ufer domains are verified in Subsection \ref{ss:pruefer}.
\smallskip

 %Inserted by TeXtelmExtel

\noindent
\textit{Proof of Theorem \ref{dimension equality}.}
Let $(O,M)\in\bL (S)$ be essentially finitely presented over $S$. Theorem \ref{upper bounds} shows that regularity of $O$ is equivalent to the finiteness of the weak dimension $\wdim O$.

 %Inserted by TeXtelmExtel

To verify the formula for the weak dimension we perform an induction starting with the case $\wdim O=1$: Lemma \ref{wdim O=1} yields $\Rad (NO)=M$, consequently  $\dim\overline{O}=0$ as asserted.

 %Inserted by TeXtelmExtel

Assume next that $\wdim O>1$ holds: the ring $\overline{O}$ is noetherian, Theorem \ref{upper bounds} thus yields $\wdim O\leq\dim\overline{O}+1$, hence $\dim\overline{O}\geq 1$. We choose a foreimage $t\in M$ of some $\overline{t}\in\overline{M}\setminus\overline{M}^2$ under the natural homomorphism $O\rightarrow\overline{O}$. We get $t\not\in M^2$, hence an application of Theorem \ref{regular prime} yields that the local ring $O/tO$ is coherent, regular, and possesses weak dimension 
\begin{equation}
\label{reduction}
\wdim (O/tO)=\wdim O -1.
\end{equation}

 %Inserted by TeXtelmExtel

\textbf{Claim:} $O/tO$ is essentially finitely presented over $S/p$, $p:=tO\cap S$, and $O/tO\in\bL (S/p)$ -- note that $S/p$ can be a field.

 %Inserted by TeXtelmExtel

\noindent
The domain $O/tO$ is a local extension of the valuation domain $S/p$ hence faithfully flat. We get $O/tO\in\bL (S/p)$. By assumption $O=A_q$, $q\in\Spec A$, for some finitely presented $S$-algebra $A$. Consequently $O/tO=B_{q_B}$ for $B:=A/(tO\cap A)$ and $q_B:=q/(tO\cap A)$. The domain $B$ is a finitely generated $S/p$-algebra and thus finitely presented -- \cite{Nag}, Thm.~3, \cite{RG}, Cor.~3.4.7.

 %Inserted by TeXtelmExtel

We are now in the position to perform the induction step in which we have to distinguish between two cases:

 %Inserted by TeXtelmExtel

\textbf{$S/p$ is no field:} by induction hypothesis we then get
\begin{equation}
\label{step}
\wdim (O/tO)=\dim (O/tO\otimes_{S/p}S/N)+1,
\end{equation}
where $O/tO\otimes_{S/p}S/N=\overline{O}/\overline{t}\,\overline{O}$, $\overline{t}:=t+NO$. The Principal Ideal Theorem yields the inequality $\dim (\overline{O}/\overline{t}\,\overline{O})\geq\dim\overline{O}-1$. The combination of the equations (\ref{reduction}) and (\ref{step}) thus leads to $\wdim O\geq\dim\overline{O}+1$. The reversed inequality $\wdim O\leq\dim\overline{O}+1$ is given by Theorem  \ref{upper bounds}.

 %Inserted by TeXtelmExtel

\textbf{$S/p$ is a field:} the ring $O/tO$ then is noetherian, therefore the homological dimensions satisfy $\wdim O/tO=\gldim O/tO$. Equation (\ref{reduction}) thus yields 
\begin{equation}
\label{almost}
\wdim O =\wdim (O/tO)+1=\gldim (O/tO)+1=\dim (O/tO)+1.
\end{equation}
The local ring $O_{tO}$ is regular, coherent and by formula (\ref{wdim O=fdim k}) has weak dimension $\wdim O=1$. Since $O_{tO}\in\bL (S)$ is essentially finitely presented over $S$, an application of Lemma \ref{wdim O=1} yields the minimality of the prime $tO$ among the primes of $O$ containing $NO$. In Subsection \ref{ss:pruefer} we shall see that the ring $\overline{O}=O/NO$ is equidimensional -- see Corollary \ref{equidimensional}. We get $\dim\overline{O}=\dim (O/tO)$ and thus  the assertion by plugging into equation (\ref{almost}).
\eproof
% ===================================
\subsection{Upper bounds for the weak dimension}
\label{ss:bounds}
% ===================================
In the present subsection we provide a proof of Theorem \ref{upper bounds} essentially by applying the following variation of a result of W.~Vasconcelos \cite{Vas}, Cor.~5.12:
\begin{proposition}
\label{wdim and radical}
Let $(O,M)$ be a coherent, regular, local ring with the property $M=\Rad I$ for some finitely generated ideal $I\lhd O$. If $I$ is generated by $\ell\in\mathbb{N}$ elements, then $\wdim O\leq\ell$ holds. If $I$ is generated by a regular sequence of length $\ell$, then the weak dimension satisfies $\wdim O=\ell$.
\end{proposition}
\proof
Let $M=\Rad I$ for some ideal $I\lhd O$ generated by $\ell\in\mathbb{N}$ elements. The equation $\wdim O=\Grade M$ (\cite{Gla2}, Lemma 3) reduces the proof of the first statement in the theorem to verifying the inequality $\Grade M\leq\ell$. According to \cite{Nor}, Ch.~5, Thm.~12 for every ideal $J\lhd O$ the equation
\begin{equation}
\label{Grade and Rad}
\Grade J=\Grade (\Rad J)
\end{equation}
holds. Moreover $\Grade J\leq\ell$ whenever $J$ is generated by $\ell$ elements 
(\cite{Nor}, Ch.~5, Thm.~13). Hence we get $\Grade M=\Grade I\leq\ell$ as asserted. If $I$ is generated by a regular sequence of length $\ell$, then we clearly have $\Grade I\geq \ell$, which yields the second assertion.
\eproof
Within the subsequent proof and at other points of this article we need to consider those primes $p\lhd S$ of a valuation domain $S$ that have the property
\begin{equation}
\label{limit prime}
p=\bigcup\limits_{p^\prime\in\Spec S:\;p^\prime\subset p}p^\prime.
\end{equation}
We call a prime $p$ with the property (\ref{limit prime}) a \textbf{limit prime of $S$}. Note that instead of taking the union of all primes $p^\prime\subset p$ we can use only those primes $p^\prime$ of the form $\Rad (tS)$ for some $t\in S$ -- in the sequel we denote this set by $\mathcal{R}(S)$.
\medskip

 %Inserted by TeXtelmExtel

\noindent
\textit{Proof of Theorem \ref{upper bounds}.} 
Let $(O,M)\in\bL ((S,N))$ be regular; we first prove the inequality $\wdim O\leq\dim\overline{T}+1$. Choose foreimages $\overline{t}_1,\ldots ,\overline{t}_\ell\in\overline{M}$ of a basis of $\overline{T}$. For each $\overline{t}_i$ choose a foreimage $t_i\in M$ under the natural homomorphism $O\rightarrow\overline{O}$. The ideal $I:=\sum\limits_{i=1}^\ell Ot_i$ then satisfies the equation $M=I+NO$ since by Nakayama's Lemma $\overline{M}$ is generated by $\overline{t}_1,\ldots ,\overline{t}_\ell$.

 %Inserted by TeXtelmExtel

If $N$ is no limit prime, then $N=\Rad (tS)$ for some $t\in N$. The ideal $J:=I+tO$ is generated by $\ell +1$ elements and satisfies $\Rad J=M$. Proposition \ref{wdim and radical} thus yields $\wdim O\leq \ell +1$ as asserted.

 %Inserted by TeXtelmExtel

In the remaining case the comment subsequent to equation (\ref{limit prime}) for the maximal ideal $N$ gives us
\[
M=\bigcup\limits_{p\in{\mathcal{R}(S)}}(I+pO).
\]
Using \cite{Gla2}, Lemma 3 and a basic property of $\Grade$ yields $\wdim O=\Grade M=\sup (\Grade J\sep J\lhd O\mbox{ finitely generated})$. Hence in order to prove the assertion in the present case it suffices to verify the inequality $\Grade J\leq\ell +1$ for every finitely generated ideal $J\lhd O$. Such an ideal is contained in some $I+pO$, $p\in\mathcal{R}(S)$. Utilizing the relation (\ref{Grade and Rad}) we get 
\[
\Grade (I+pO)\leq\Grade (\Rad (I+tO))=\Grade (I+tO)\leq\Grade (I+pO) 
\]
and therefore $\Grade J\leq\Grade (I+pO)=\Grade (I+tO)\leq\ell +1$.

 %Inserted by TeXtelmExtel

In the case of a noetherian ring $\overline{O}$ we start by choosing a system of parameters  $\overline{t}_1,\ldots ,\overline{t}_\ell\in\overline{M}$, $\ell=\dim\overline{O}$. The ideal $I:=\sum\limits_{i=1}^\ell Ot_i$ generated by foreimages $t_i\in M$ of the elements $\overline{t}_i$ under the map $O\rightarrow\overline{O}$ then satisfies $\Rad (I+NO)=M$ since by construction $\Rad ((I+NO)/NO)=\overline{M}$. From here on the proof  proceeds almost identically to the proof of the preceeding case.
\eproof
% =========================================
\subsection{Properties of algebras over Pr\"ufer domains}
\label{ss:pruefer}
% =========================================
Let $R$ be a Pr\"ufer domain with field of fractions $\Frac R=:K\neq R$. In the sequel we discuss various properties of $R$-algebras, that are relevant for the present article. Partially they have already been used in Subsection \ref{ss:main}.

 %Inserted by TeXtelmExtel

Every finitely generated ideal $I\lhd R$ is invertible and thus finitely presented -- a Pr\"ufer domain is coherent. While in general a finitely presented algebra over a coherent ring needs not be coherent, Pr\"ufer domains show a smoother behavior in that respect:
\begin{theorem}
\label{coherence}
Every finitely generated, flat algebra $A$ over a Pr\"ufer domain $R$ is finitely presented and coherent.
\end{theorem}
\proof
A finitely generated, flat algebra over a domain is finitely presented -- \cite{RG}, Cor.~3.4.7. The coherence of $A$ now follows from the coherence of the polynomial ring $R[X_1,..,X_n]$ over a Pr\"ufer ring (\cite{Sab}, Prop.~3, \cite{Gla1}, Cor.~7.3.4.) and the coherence of factor rings of a coherent ring with respect to a finitely generated ideal.
\eproof
\begin{corollary}
Let $(O,M)$ be a local, flat, essentially finitely presented $R$-algebra and $p:=M\cap R$, then $O\in\bL (R_p)$.
\end{corollary}
\proof
Choose a finitely presented $R$-algebra $A$ such that $O=A_q$, $q\in\Spec A$, holds. The set of primes $q^\prime\in\Spec A$ such that $A_{q^\prime}|R_{q^\prime\cap R}$ is flat, is open and by assumption non-empty. Thus we can replace $A$ by a finitely presented, flat $R$-algebra. The stability of coherence under localization and Theorem \ref{coherence} yield the coherence of $O$ and thus $O\in\bL(R_p)$.
\eproof
A coherent, regular, local ring $O$ is a domain -- \cite{Gla1}, Lemma 4.2.3. The subsequent auxiliary result is motivated by this fact and has been applied in the proof of Theorem \ref{dimension equality}.
\begin{lemma}
\label{domain}
Let $R$ be a non-trivial valuation domain. Every essentially finitely presented domain $O\in\bL (R)$ can be written as $O=A_q$ for some finitely generated domain $A$.
\end{lemma}
\proof
Choose a finitely presented, flat $R$-algebra $A$ such that $O=A_q$ for some $q\in\Spec A$. The flatness of $A|R$ implies $q_0\cap R=0$ for every minimal prime ideal $q_0\subseteq q$. Since $O$ is a domain we have $A_q=(A/q_0)_{q/q_0}$.
\eproof
We finish this subsection with a dimension-theoretic property of algebras over Pr\"ufer domains:
\begin{theorem}
\label{affine fibre dimension}
Let $A$ be a domain that is finitely generated over the Pr\"ufer domain $R$. Then the rings $A\otimes_R kp$, $p\in\Spec R$, $kp=\Frac (R/p)$, are equidimensional of dimension $\dim (A\otimes_R K)$, $K=\Frac R$.
\end{theorem}
\proof
Let $q\in\Spec A$ be minimal among the primes lying over $p\in\Spec R\setminus 0$. \cite{Nag}, Lemma 2.1 yields 
$\trdeg (A|R)=\trdeg (A/q|R/p)$. The assertion now follows from the fact that $\trdeg (A|R)=\dim (A\otimes_R K)$ and $\trdeg (A/q|R/p)=\dim (A/q\otimes_{R/p} kp)$ hold.
\eproof
\begin{corollary}
\label{equidimensional}
Let $(S,N)$ be a non-trivial valuation domain. The factor ring $\overline{O}=O/NO$ of an essentially finitely presented, local domain $O\in\bL ((S,N))$ is equidimensional.
\end{corollary}
\proof
Lemma \ref{domain} shows that $O=A_q$, $q\in\Spec A$, for some domain $A$ finitely generated over $S$. 
\eproof
% ===================
\section{Regular sequences}
\label{s:regular sequences}
% ===================
In the present section we are concerned with the construction of maximal regular sequences within the maximal ideal $M$ of a coherent local ring $O$ of finite weak dimension. These sequences are obtained through lifting linearly independent elements from the $O/M$-vector space $M/M^2$ to $M$. We study this approach in Subsection \ref{ss:regseq+wdim}, where we also provide the proof of Theorem \ref{regular prime} that was postponed in Subsection \ref{ss:main}. In Subsection \ref{ss:maxregseq} we apply the results to a regular, local, essentially finitely presented, flat $S$-algebra $O$ over the valuation domain $S$. We shall see that in this case regular sequences of the maximal possible length $\wdim O$ exist within $M$. As a consequence the fibre ring $O/(M\cap S)O$ is shown to be Cohen-Macaulay, and even regular in the case that $M\cap S$ is {\em not} finitely generated.
% =========================================
\subsection{Coherent local rings of finite weak dimension}
\label{ss:regseq+wdim}
% =========================================
For a local ring $(O,M)$ the $O$-module $T:=M/M^2$ forms a vector space over the residue field $k:=O/M$. If $O$ is a noetherian regular ring, then a regular sequence $(t_1,\ldots ,t_\ell)$ in $M$ can be constructed by choosing the elements $t_i\in M$, $i=1,\ldots ,\ell$, to be $k$-linearly independent modulo $M^2$. The ideal $q\lhd O$ generated by $t_1,\ldots ,t_\ell$ is prime and the factor ring $O/q$ is regular too. In the sequel we show that this method carries over to coherent local rings of finite weak dimension:
\begin{theorem}
\label{regular sequence}
In a coherent local ring $(O,M)$ of finite weak dimension every set $\{t_1,\ldots ,t_\ell\}\subset M$ such that $\{t_1+M^2,\ldots ,t_\ell +M^2\}\subset T$ is $k$-linearly independent has the properties:
\begin{enumerate}
\item $(t_1,\ldots ,t_\ell)$ is a regular sequence,
\item $q:=\sum\limits_{i=1}^\ell Ot_i$ is a prime ideal,
\item $\wdim O_q=\ell$ and $\wdim O/q=\wdim O-\ell$.
\end{enumerate}
\end{theorem}
\noindent
Note however that in the situation given in Theorem \ref{regular sequence} the vector space $T$ may have strange properties:
\begin{example}
\label{strange T}
If the maximal ideal of a valuation domain $(S,N)$ is not finitely generated it satisfies $N=N^2$, thus $T=0$ holds although $\grade N=1$. The localization $O:=S[X]_q$ of the polynomial ring $S[X]$ at $q:=XS[X]+N[X]$ is coherent (\cite{Gla1}, Thm.~7.3.3) and satisfies $\wdim O=2$ according to Hilbert's Szyzygy Theorem \cite{Vas}, Thm.~0.14. On the other hand $T=k(X+M^2)$.
\end{example}%
This example demonstrates that we cannot expect to obtain maximal regular sequences within $M$ through lifting a $k$-basis of $T$ as in the noetherian case. Nevertheless we pursue the approach of lifting as far as possible. 

 %Inserted by TeXtelmExtel

As one expects the proof of Theorem \ref{regular sequence} is an induction on the number $\ell$ of $k$-linearly independent elements. The case of one element $t\in M\setminus M^2$ is treated in Theorem \ref{regular prime} whose proof was postponed in Subsection \ref{ss:main}; we provide this proof now:
\smallskip

 %Inserted by TeXtelmExtel

\noindent
\textit{Proof of Theorem \ref{regular prime}.}
Every coherent, regular, local ring is a domain \cite{Gla1}, Cor.~6.2.4, a fact that we need twice throughout the current proof.

 %Inserted by TeXtelmExtel

The weak dimension of the coherent ring $O/tO$ can be computed through formula (\ref{wdim O=fdim k}) in Subsection \ref{ss:regularity}: $\wdim (O/tO)=\fdim_{O/tO} (k)$, $k=O/M$. In order to prove regularity of $O/tO$ it therefore suffices to show $\fdim_{O/tO} (M/tO)\neq\infty$.

 %Inserted by TeXtelmExtel

\textbf{Claim:} $\fdim_{O/tO} (M/tM)\neq\infty$ implies $\fdim_{O/tO} (M/tO)\neq\infty$.

 %Inserted by TeXtelmExtel

Choose a set $\mathcal{B}\subset T$ such that $\mathcal{B}\cup\{\overline{t}\}, \overline{t}:=t+M^2$ is a $k$-basis of $T$. For each $\overline{b}\in\mathcal{B}$ choose a foreimage in $M$; let $B\subset M$ be the set of these foreimages. The ideal
\[
I:=\sum\limits_{b\in B}bO +M^2
\]
has the properties $M=I+tO$ and $I\cap tO=tM$. Only the inclusion $\subseteq$ of the second property requires a verification: an element $a\in I\cap tO$ can be expressed as
\[
tr=a=\sum\limits_{b\in B}r_b b+s,\; r,r_b\in O,\; s\in M^2.
\]
Taking this equation modulo $M^2$ and using the $k$-linear independence of $\mathcal{B}\cup\{\overline{t}\}$ yields $r,r_b\in M$ as desired.

 %Inserted by TeXtelmExtel

Utilizing the properties of $I$ we see that the exact sequence
\[
0\rightarrow tO/tM\rightarrow M/tM\rightarrow M/tO\rightarrow 0
\]
splits and we obtain $M/tM\cong M/tO\oplus tO/tM$ as $O/tO$-modules. Since for arbitrary modules $N,P$ the inequality $\fdim N\leq\fdim (N\oplus P)$ holds, the claim is proved.

 %Inserted by TeXtelmExtel

We now have to verify $\fdim_{O/tO} (M/tM)\neq\infty$: the maximal ideal $M$ by assumption possesses a flat resolution
\[
0\rightarrow F_m\rightarrow\ldots\rightarrow F_0\rightarrow M\rightarrow 0
\]
of finite length. Tensoring with $O/tO$ yields a resolution
\[
0\rightarrow F_m\otimes_O (O/tO)\rightarrow\ldots\rightarrow F_0\otimes_O (O/tO)\rightarrow M/tM\rightarrow 0
\]
of $M/tM$ by flat $O/tO$-modules; here we use the fact that $t$ is no zero-divisor.

 %Inserted by TeXtelmExtel

The formula for the weak dimension of $O/tO$ follows from \cite{Gla1}, Thm.~3.1.4 (2), where we again need the fact that $t$ is no zero-divisor.
\eproof
\textit{Proof of Theorem \ref{regular sequence}.}
We perform an induction on the number $\ell$ of elements in the sequence.

 %Inserted by TeXtelmExtel

In the case $\ell=1$, point 1 is clear and point 2 as well as the second part of point 3 are the content of Theorem \ref{regular prime}. An application of formula (\ref{wdim O=fdim k}) in Subsection \ref{ss:regularity} yields $\wdim O_{t_1O}=1$, since the maximal ideal of $O_{t_1O}$ is principal.

 %Inserted by TeXtelmExtel

Assume now that $\ell >1$ holds and consider the coherent, regular, local ring $O^\prime :=O/t_1O$ of weak dimension $\wdim O^\prime =\wdim O-1$.\\ The elements $t_i^\prime :=t_i+t_1O$, $i=2,\ldots ,\ell$, are $k$-linearly independent modulo $(M^\prime )^2$, where $M^\prime :=M/t_1O$. Therefore by induction hypothesis the sequence  $(t_2^\prime ,\ldots ,t_\ell^\prime )$ is regular, so that the sequence $(t_1,\ldots ,t_\ell)$ is regular too. Moreover the ideal $q:=\sum\limits_{i=1}^\ell Ot_i$ is prime since this holds for $q^\prime :=q/t_1O\lhd O^\prime$; we get
\[
\wdim (O/q)=\wdim (O^\prime /q^\prime )=\wdim O^\prime -(\ell -1)=\wdim O -\ell .
\]
As a localization of a regular ring $O^\prime_{q^\prime}=O_q/t_1O_q$ is regular too. We apply \cite{Gla1}, Thm.~3.1.4 (2) and obtain
\[
\ell -1=\wdim O^\prime_{q^\prime}=\wdim O_q -1
\]
as asserted.
\eproof
\begin{corollary}
\label{inequalities}
For every coherent local ring $(O,M)$ of finite weak dimension the inequalities
\[
\dim_k T\leq\grade M\leq\wdim O
\]
hold.
\end{corollary}
\begin{remark}
We have already seen in Example \ref{strange T} that the first inequality in Corollary \ref{inequalities} can be strict. An example for the strict inequality $\grade M<\wdim O$ seems to be missing in the literature. 
\end{remark}
\begin{corollary}
\label{equality}
For a coherent local ring $(O,M)$ of finite weak dimension the equation $\dim_k T=\wdim O
$ holds if and only if $M$ is finitely generated. $M$ is then generated by a regular sequence of length $\wdim O$.
\end{corollary}
\proof
Let $q\lhd O$ be an ideal generated by a set $t_1,\ldots ,t_\ell$ of foreimages of a $k$-basis of $T$.
Assume that $\dim_k T=\wdim O$ holds. Theorem \ref{regular sequence} then yields $\wdim O/q=0$, which implies that $O/q$ is a field. Assume now that $M$ is finitely generated. The regular local ring $O^\prime:=O/q$ has the property $(M/q)/(M/q)^2=0$. Since $M/q$ is finitely generated Nakayama's Lemma yields $q=M$.
\eproof
\begin{remark}
\label{chinese}
Based on the work of Vasconcelos \cite{Vas} the regularity of a coherent local ring $(O,M)$ with finitely generated maximal ideal $M$ has been investigated in \cite{TZT}. One of the main results of this article states that the regularity of $O$ is equivalent to the existence of a regular sequence generating $M$.  In \cite{THT} this equivalence is generalized to indecomposable, semilocal, coherent rings.
\end{remark}
% ===========================================
\subsection{Regular local algebras over a valuation domain}
\label{ss:maxregseq}
% ===========================================
In a valuation domain $(S,N)$ radical ideals are prime and the primes of the form $p=\Rad (tS)$, $t\in S$, are precisely the non-limit primes -- see (\ref{limit prime}) in Subsection \ref{ss:bounds}. Utilizing this property of prime ideals Theorem \ref{regular sequence} applied to a regular, essentially finitely presented algebra $(O,M)\in\bL (S)$ yields distinguished regular sequences within $M$:
\begin{theorem}
\label{maximal regular sequences}
Let $(S,N)$ be a non-trivial valuation domain. The maximal ideal of a regular, essentially finitely presented $S$-algebra $(O,M)\in\bL (S)$ contains a maximal regular sequence $(t_1,\ldots ,t_d)$ of length $d=\wdim O$. 

 %Inserted by TeXtelmExtel

More precisely if $t_1,\ldots ,t_\ell\in M$ are elements such that  $(t_1+M^2,\ldots ,t_\ell+M^2)$ forms a $k$-basis of $T=M/M^2$, then:
\begin{enumerate}
\item If $M$ is finitely generated, then $(t_1,\ldots ,t_\ell )$ is a maximal regular sequence of length $\ell=\wdim O$ and $M=\sum\limits_{i=1}^\ell Ot_i$.
\item If $M$ is not finitely generated,  then for every $t\in N\setminus 0$ the sequence $(t_1,\ldots ,t_\ell,t)$ is maximal regular of length $\ell +1=\wdim O$. The ideal $\Rad (\sum\limits_{i=1}^\ell Ot_i +Ot)$ is a prime ideal.
If $t\in N$ satisfies $\Rad (tS)=N$, then the equation $M=\Rad (\sum\limits_{i=1}^\ell Ot_i +Ot)$ holds.
\end{enumerate}
\end{theorem}
\begin{remark}
\label{length maximal regular sequence}
Note that Theorem \ref{maximal regular sequences} implies $\grade M=\wdim O$ for regular, essentially finitely presented $S$-algebras $(O,M)\in\bL (S)$. In \cite{Vas} Vasconcelos showed that a maximal regular sequence within a coherent, regular, local ring $(O,M)$ satisfying $\grade M=\wdim O$ may have a length smaller than $\wdim O$. The author does not know whether this can happen for the local rings considered in the theorem.
\end{remark}%
As an important consequence of Theorem \ref{maximal regular sequences} we obtain information about the structure of the fibre ring $\overline{O}=O/NO$ and get a result that is analogous to the noetherian case if $N$ is finitely generated but possesses a surprising non-noetherian component otherwise:
\newpage
\begin{theorem}
\label{fibre ring}
The fibre ring $\overline{O}$ of a regular, local, essentially finitely presented $S$-algebra $(O,M)\in\bL (S)$ over a non-trivial valuation domain $(S,N)$ has the properties:
\begin{enumerate}
\item $\overline{O}$ is Cohen-Macaulay.
\item If $N$ is not finitely generated, then $\overline{O}$ is regular.
\end{enumerate}
\end{theorem}
\proof
For a non-finitely generated ideal $N\lhd S$ the second assertion implies the first. We thus prove the first assertion only in the case $N=tS$.

 %Inserted by TeXtelmExtel

We have to show that $\dim\overline{O}=\grade\overline{M}$ holds; to this end we apply the formula \cite{Alf}, Cor.~2.10 that describes the behavior of polynomial grade in essentially finitely presented ring extensions like $O|S$:
\begin{equation}
\label{Grade in extension}
\Grade_O M=\Grade_S N+\Grade_{\,\overline{O}}\overline{M};
\end{equation}
the indices denote the base ring with respect to which polynomial grade has to be calculated. Theorem \ref{dimension equality} yields $\Grade_O M=\wdim O=\dim\overline{O}+1$. Furthermore $\Grade_S N=1$ since $S$ is a valuation domain. Finally $\Grade_{\,\overline{O}}\overline{M}=\grade_{\,\overline{O}}\overline{M}$ because $\overline{O}$ is noetherian. Equation (\ref{Grade in extension}) consequently leads to the claim $\dim\overline{O}=\grade_{\,\overline{O}}\overline{M}$.

 %Inserted by TeXtelmExtel

Assume next that $N$ is not finitely generated and consider the natural homomorphism
\begin{equation}
\label{cotangent reduction}
\tau :\;T\rightarrow\overline{T}:=(M/NO)/(M/NO)^2.
\end{equation}
Since by assumption $N=N^2$ holds we get $\Ker\tau =(M^2+NO)/M^2=0$. The surjectivity of $\tau$ thus yields $\dim_k T=\dim_k\overline{T}$, $k=O/M$. On the other hand by Theorem \ref{maximal regular sequences}, (2) and Theorem \ref{dimension equality}
\[
\dim_k T=\wdim O-1=\dim (\overline{O})
\]
and thus $\dim (\overline{O})=\dim_k\overline{T}$, which implies the regularity of $\overline{O}$.
\eproof
\begin{remark}
Alfonsi gives a homological definition of non-noetherian grade for arbitrary modules. The proof for the equivalence of his definition with the one given by Northcott in the case of ideals can be found in \cite{Gla1}, Thm.~7.1.8.
\end{remark}%
The remainder of the subsection is devoted to the proof of Theorem \ref{maximal regular sequences}.

 %Inserted by TeXtelmExtel

Consider a set $P:=\{t_1,\ldots ,t_\ell\}\subset M$ such that $t_1+M^2,\ldots ,t_\ell +M^2\in T$ are $k$-linearly independent. According to Theorem \ref{regular sequence} the sequence $(t_1,\ldots ,t_\ell )$ is regular and the ideal
\begin{equation}
\label{def regular prime}
q_P:=\sum\limits_{i=1}^\ell Ot_i
\end{equation}
is a prime ideal of $O$. The proof of Theorem \ref{maximal regular sequences} strongly relies on the particular distribution of these primes within $\Spec O$.
\begin{proposition}
\label{distribution regular primes}
Let $(S,N)$ be a non-trivial valuation domain and assume that $(O,M)\in\bL (S)$ is essentially finitely presented and regular. Then the prime ideals defined by equation (\ref{def regular prime}) share the property $q_P\cap S\in\{0,N\}$. The equation $q_P\cap S=N$ implies that $N$ is principal.
\end{proposition}
\proof
We perform an induction on the cardinality $\ell$ of $P$. 

 %Inserted by TeXtelmExtel

Assume first that $q_P=tO$ and let $p:=tO\cap S$; moreover assume that $p\neq 0$ holds. The regular local ring $O_{tO}\in\bL (S_p)$ is essentially finitely presented over $S_p$. Equation (\ref{wdim O=fdim k}) in Subsection \ref{ss:regularity} yields $\wdim O_{tO}=1$, that is $O_{tO}$ is a valuation domain. 
\smallskip

 %Inserted by TeXtelmExtel

\textbf{Claim}: there exists a prime $Q\in\Spec O$ containing $tO$ with the properties $Q\cap S=N$ and $O_Q$ is a valuation domain.

 %Inserted by TeXtelmExtel

We choose a domain $A$ finitely generated over $S$ such that $O=A_q$ for some prime $q\lhd A$ (Lemma \ref{domain}). Lemma \ref{wdim O=1} shows that the prime $q_t:=tO\cap A$ is minimal among the primes of $A$ lying over $p$. An application of \cite{Nag}, Lemma 2.1 yields:
\begin{equation}
\label{trdeg1}
\trdeg (A/q_t |S/p)=\trdeg (A|S).
\end{equation}
Let $\overline{q}_0\in\Spec (A/q_t)$ be minimal among the primes containing $N/p$; $\overline{q}_0$ can be choosen such that $\overline{q}_0\subseteq q/q_t$ holds. Again utilizing \cite{Nag}, Lemma 2.1 the prime $q_0\in\Spec A$ with $q_0/q_t=\overline{q}_0$ satisfies:
\begin{equation}
\label{trdeg2}
\trdeg (A/q_0 |S/N)=\trdeg (A/q_t |S/p).
\end{equation}
Combining (\ref{trdeg1}) and (\ref{trdeg2}) we can choose a transcendence basis $(x_1,\ldots ,x_m)$ of $\Frac A|\Frac S$ within $A$ such that  $(x_1+q_0,\ldots ,x_m+q_0)$ forms a transcendence basis of $\Frac A/q_0 |\Frac S/N$. This choice yields
\[
q_0\cap S[x_1,\ldots ,x_m]=N[x_1,\ldots ,x_m].
\]
The localization $B:=S[x_1,\ldots ,x_m]_{N[x_1,\ldots ,x_m]}$ is a valuation domain contained in the local ring $A_{q_0}=O_Q$, $Q:=q_0O$. Coherent regular rings are normal (\cite{Gla1}, page 205), thus the coherent regular ring $O_Q$ contains the integral closure $B^\prime$ of $B$ in $\Frac A$. Since $B^\prime$ is a Pr\"ufer domain $O_Q$ is indeed a valuation domain.

 %Inserted by TeXtelmExtel

In a valuation domain a finitely generated prime either equals $0$ or the maximal ideal. Consequently we get $tO_Q=QO_Q$ and thus $tO\cap S=N$. As we have seen in the proof of the preceeding claim the valuation domain $O_Q$ dominates the valuation domain $B$ and the extension $\Frac O_Q|\Frac B$ is finite. General valuation theory yields that finite generation of the maximal ideal $QO_Q$ descends to the maximal ideal $N[x_1,\ldots ,x_m]B\lhd B$. Let $f\in N[x_1,\ldots ,x_m]$ be a generator of $N[x_1,\ldots ,x_m]B$ and let $I\lhd S$ be generated by the coefficients of $f$, then $I=cS$ for some $c\in I$ and the factorization $f=sf^\ast$, where $f^\ast\not\in N[x_1,\ldots ,x_m]$ eventually yields $N=I$ as asserted.

 %Inserted by TeXtelmExtel

We now turn to the case $P=\{t_1,\ldots ,t_\ell\}$ with $\ell >1$.

 %Inserted by TeXtelmExtel

We already know that $t_1O\cap S\in\{0,N\}$. Since the equation $t_1O\cap S=N$ forces $q_P\cap S=N$ and since in this case we also know that $N$ must be principal  it suffices to treat the case $t_1O\cap S=0$.

 %Inserted by TeXtelmExtel

Consider the set $P_1:=\{t_2+t_1O,\ldots ,t_\ell +t_1O\}\subset O/t_1O$, and recall that $O/t_1O\in\bL (S)$ is regular. The elements of $P_1$ are $k$-linearly independent modulo $(M/t_1O)^2$. By induction hypothesis we can thus assume that the prime ideal 
\[
q_{P_1}:=\sum\limits_{i=2}^\ell (O/t_1O)(t_i+t_1O)\lhd O/t_1O
\]
satisfies $q_{P_1}\cap S\in\{0,N\}$ and that in the second case $N$ is principal. The assertions about $q_P$ now follow.
\eproof
\begin{remark}
The conclusions of Proposition \ref{distribution regular primes} hold for every principal prime ideal $tO$. The requirement $t\not\in M^2$ is not needed in the proof of that particular case.
\end{remark}%
\textit{Proof of Theorem \ref{maximal regular sequences}.}
Let $P:=\{ t_1,\ldots ,t_\ell \}$ be choosen such that its elements modulo $M^2$ form a $k$-basis of $T$.

 %Inserted by TeXtelmExtel

If $M$ is finitely generated, then $M=q_P$ holds by Nakayama's lemma. Moreover by Theorem \ref{regular sequence}, (3), the sequence $(t_1,\ldots ,t_\ell )$ possesses the maximal possible length $\ell =\wdim O$.

 %Inserted by TeXtelmExtel

Assume now that $M$ and thus $N$ are not finitely generated. We have seen in the proof of Theorem \ref{fibre ring} that in this case the natural map $\tau :T\rightarrow\overline{T}$ is an isomorphism -- see (\ref{cotangent reduction}). As a consequence the elements $t_1+NO,\ldots ,t_\ell +NO$ generate $\overline{M}$, hence $M=q_P+NO$ holds. We next first apply Theorem \ref{regular sequence} followed by Theorem \ref{upper bounds} to obtain 
\[
\wdim (O/q_P)\leq\dim\overline{T}+1-\dim T.
\]
Hence using $\dim T=\dim\overline{T}$ we arrive at the inequality $\wdim (O/q_P)\leq 1$. The equation $\wdim (O/q_P)=0$ implies that $O/q_P$ is a field and thus that $M=q_P$ is finitely generated, which has been excluded. Consequently $O/q_P$ is a valuation domain $S^\prime$ and $\ell=\wdim O-1$ by Theorem \ref{regular sequence} (3).

 %Inserted by TeXtelmExtel

Proposition \ref{distribution regular primes} tells us that $q_P\cap S=0$, because $N$ is not finitely generated. We conclude that $S^\prime$ is a local extension of $S$ and that therefore $(t_1,\ldots ,t_\ell,t)$ is a regular sequence for every non-zero $t\in N$. Moreover this sequence possesses the maximal possible length $\wdim O$.

 %Inserted by TeXtelmExtel

If $N=\Rad (tS)$ for some some $t\in N$, then $M=q_P+NO=\Rad (q_P+tO)$ as asserted. For an arbitrary $t\in N\setminus 0$ the ideal $\Rad (tS^\prime )$ is prime, hence the same holds for $\Rad (q_P+tO)$.
\eproof
\begin{corollary}
\label{descend of finite generation}
Under the assumptions made in Theorem \ref{maximal regular sequences} for the local extension $(O,M)|(S,N)$ the maximal ideal $M$ is finitely generated if and only if $N$ is finitely generated (that is principal).
\end{corollary}
\newpage
% ======================================
\section{Families over Pr\"ufer domains}
\label{s:geometry}
% ======================================
The results presented so far in this article grew out of an attempt to understand the geometric structure of certain non-noetherian schemes:
\begin{definition}
\label{R-family}
Let $R$ be a Pr\"ufer domain. An \textbf{$R$-family} is a separated, faithfully flat $R$-scheme $\cX$ of finite type such that the fibres $\cX_p:=\cX\times_R kp$, $p\in\Spec R$, are equidimensional $kp$-schemes of a common dimension $n\in\mathbb{N}$. An {\em integral} $R$-family $\cX$ is called \textbf{$R$-variety}. The common dimension $n$ of the fibres is called the \textbf{relative dimension of $\cX$ over $R$}.
\end{definition}%
Throughout this section assume again $R$ to be a Pr\"ufer domain with field of fractions $K\neq R$. As an immediate consequence of the Theorems \ref{coherence} and \ref{affine fibre dimension} we get:
\begin{theorem}
\label{finite type}
Let $R$ be a Pr\"ufer domain.
\begin{enumerate}
\item An $R$-family $\cX$ is of finite presentation and its structure sheaf $\cO_{\cX}$ is coherent.
\item A separated faithfully flat $R$-scheme $\cX$ of finite type possessing an equidimensional generic fibre is an $R$-family.
\end{enumerate}
\end{theorem}
For an arbitrary scheme $\cX$ we define the \textbf{regular locus $\Reg\cX$} to be the set of points $P\in\cX$ such that the local ring $\cO_{\cX,P}$ at $P$ is regular (in the sense of Bertin). Theorem \ref{dimension equality} allows to give a description of the regular locus of an $R$-family:
\begin{theorem}
\label{regular locus}
Let $\cX$ be an $R$-family of relative dimension $n\in\mathbb{N}$ over the Pr\"ufer domain $R$. The regular locus of $\cX$ is characterized through
\[
\Reg\cX =\{P\in\cX | \wdim\cO_{\cX,P}\leq n+1\}.
\]
The upper bound $n+1$ is attained at the closed points of the fibres $\cX_p$, where $p\in\Spec R\setminus 0$.
\end{theorem}
\smallskip

 %Inserted by TeXtelmExtel

\noindent
A rough estimate of the extend of the singular locus $\cX\setminus\Reg\cX$ of an $R$-family $\cX$ is provided by the next result.
\begin{theorem}
\label{regular locus and fibres}
The regular locus of an $R$-family $\cX$ satisfies:
\[
\bigcup\limits_{p\in\Spec R}\Reg\cX_p\subseteq\Reg\cX .
\]
\end{theorem}
\proof
For a point $P\in\Reg\cX_p$ the exact sequence
\[
0\rightarrow p\cO_{\cX,P}\rightarrow\cO_{\cX,P}\rightarrow\cO_{\cX_p,P}\rightarrow 0
\]
is a flat resolution of the $\cO_{\cX,P}$-module $\cO_{\cX_p,P}$: since $\cO_{\cX,P}$ is a flat $R$-module the ideal $p\cO_{\cX,P}\cong p\otimes_R\cO_{\cX,P}$ is a flat $\cO_{\cX,P}$-module. We can thus apply \cite{Gla1}, Thm.~3.1.1 and obtain:
\[
\wdim\cO_{\cX,P}\leq\wdim\cO_{\cX_p,P}+\fdim\cO_{\cX_p,P}\leq\dim\cO_{\cX_p,P}+1,
\]
where in the last inequality we use the equation $\wdim\cO_{\cX_p,P}=\dim\cO_{\cX_p,P}$ for the noetherian, regular, local ring $\cO_{\cX_p,P}$.
\eproof
It is well-known that the regular locus of a reduced scheme $Y$ of finite type over a field $k$ is open and non-empty. Therefore the significance of Theorem \ref{regular locus and fibres} is increased through the following fact proved in \cite{Kn2}: denote the set of maximal ideals of the Pr\"ufer domain $R$ by $\MaxSpec R$ and define
\[
\fgSpec R:=\{p\in\Spec R\sep pR_p\mbox{ is finitely generated.}\}.
\]
According to \cite{Kn2}, Cor.~2.8 a normal $R$-variety $\cX$ of relative dimension $1$ has the property: the fibres $\cX\times_R kp$, $p\not\in\fgSpec R$, are reduced. If $\cX|\Spec R$ is proper, then all fibres over primes $p\not\in\MaxSpec R$ are reduced. The valuation-theoretic proof given in \cite{Kn2} for this fact easily carries over to the case of arbitrary relative dimension.

 %Inserted by TeXtelmExtel

We have seen in Subsection \ref{ss:maxregseq} that the set $\fgSpec R$ plays a particular role in the context of regularity too:
\begin{theorem}
\label{fibres}
The fibres of an $R$-family $\cX$ over the Pr\"ufer domain $R$ possess the following properties:
\begin{enumerate}
\item If $\cX$ is regular, then all fibres $\cX_p$ are Cohen-Macaulay, while the fibres $\cX_p$, $p\not\in\fgSpec R$, are regular $kp$-schemes.
\item If the structure morphism $\cX\rightarrow\Spec R$ is closed and all fibres $\cX_p$, $p\in\MaxSpec R$, are regular, then $\cX$ is regular.
\end{enumerate}
\end{theorem}
\proof
Assertion 1 follows from Theorem \ref{fibre ring}.

 %Inserted by TeXtelmExtel

Assertion 2: by assumption about the structure morphism the Zariski closure of every point $P\in\cX$ intersects some fibre $\cX_p$, $p\in\MaxSpec R$. Therefore
\[
\cO_{\cX,P}=(\cO_{\cX,Q})_q,\; q\in\Spec\cO_{\cX,Q},
\]
for some $Q\in\cX_p$. By assumption $\cO_{\cX,Q}$ is regular, therefore $\cO_{\cX,P}$ is regular too.
\eproof
We conclude this section with the discussion of an example:
\begin{example}
\label{arithmetic elliptic}
Let $R$ be a $2$-dimensional valuation domain with field of fractions $K$ and assume that the spectrum of $R$ satisfies $\Spec R=\fgSpec R=\{0,p,N\}$. Choose elements $s,t\in N$ such that $pR_p=sR_p$ and $N=tR$ hold.

 %Inserted by TeXtelmExtel

Consider the projective $R$-family
\begin{equation}
\label{example}
\cC :=\Proj (R[X,Y,Z]/FR[X,Y,Z]), \;\; F:=sX^3+tZ^3-ZY^2.
\end{equation}
\begin{figure}[h]
\begin{center}
\epsfig{file=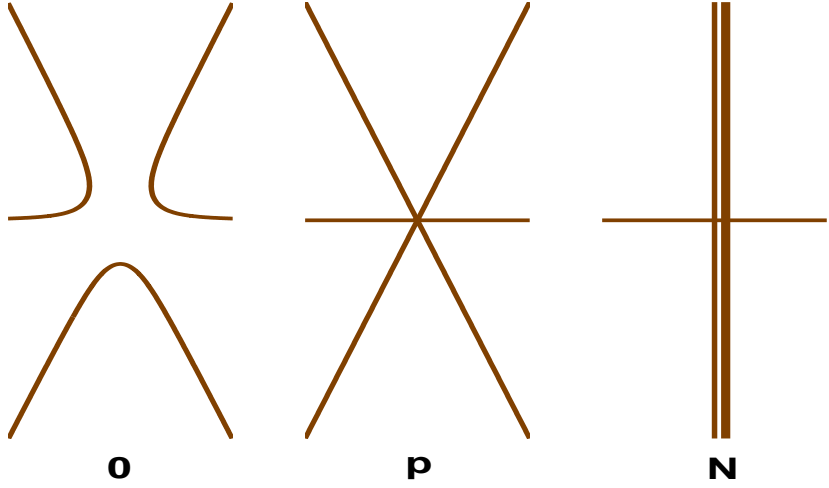,width=11cm,height=6.5cm,angle=0}
\caption{Fibres of the $R$-curve $U$}
\end{center}
\end{figure}
The generic fibre $C:=\cC\times_RK$ is an elliptic curve -- of course we have to assume $\ch K\neq 2,3$ at this point. Consequently the scheme $\cC$ is an $R$-variety of relative dimension 1 (Theorem \ref{finite type}) and $C\subset\Reg\cC$.

 %Inserted by TeXtelmExtel

The fibre $\cC_p=\cC\times_R kp$ is defined through the homogenous polynomial
\[
(t+p)Z^3-ZY^2\in kp[X,Y,Z].
\]
Assuming that $kp$ is perfect and using the Jacobian criterion we see that the point $Q:=[1:0:0]\in\mathbb{P}^2_{kp}$ is the only singular point of $\cC_p$. It lies in the affine open subscheme $U\subset\cC\subset\mathbb{P}^2_R$ defined through the condition $X\neq 0$. The coordinate ring of $U$ is given through
\begin{equation}
\label{X-chart}
R[y,z]:=R[Y,Z]/GR[Y,Z],\; G:=s+tZ^3-ZY^2\in R[Y,Z],
\end{equation}
hence
\[
\cO_{\cC,Q}=R[y,z]_q,\; q:=R[y,z]y+R[y,z]z+pR[y,z].
\]
The relation $s=zy^2-tz^3$ (\ref{X-chart}) shows that $qR_p[y,z]$ and thus the maximal ideal of $\cO_{\cC,Q}$ are generated by $y,z$.

 %Inserted by TeXtelmExtel

We claim that $(y,z)$ is a regular sequence in $R_p[y,z]$ and thus in $\cO_{\cC,Q}$: indeed
\[
R_p[y,z]/yR_p[y,z]\cong R_p[Z]/(tZ^3+s)R_p[Z]
\]
and the polynomial $tZ^3+s$ is prime in $R_p[Z]$. Consequently $y\in R_p[y,z]$ is prime and $z$ is no zero-divisor on $R_p[y,z]/yR_p[y,z]$. 

 %Inserted by TeXtelmExtel

We conclude that $Q\in\Reg\cC$: the Koszul complex of the maximal ideal $\cM_{\cC,Q}$ yields a flat resolution of length $2$, therefore $\wdim\cO_{\cC,Q}\leq 2$ by Formula (\ref{wdim O=fdim k}) in Subsection \ref{ss:regularity}.

 %Inserted by TeXtelmExtel

Since $\Reg\cC_p=\cC_p\setminus Q$ Theorem \ref{regular locus and fibres} implies $\cC_p\subset\Reg\cC$.

 %Inserted by TeXtelmExtel

The closed fibre $\cC_N:=\cC\times_R kN$ is defined through the homogenous polynomial
\[
ZY^2\in kN[X,Y,Z].
\]
The singular locus of $\cC_N\subset\mathbb{P}^2_{kN}$ thus equals the doubled line given through $Y=0$. For the sake of simplicity we assume $kN$ to be algebraically closed. The singular points of the form $[\alpha :0:1]\in\cC_N$, $\alpha\in kN$, lie in the affine open subscheme $V\subset\cC$ defined through $Z\neq 0$. The coordinate ring of $V$ thus is given through
\begin{equation}
\label{Z-chart}
R[x,y]:=R[X,Y]/HR[X,Y],\; H:=sX^3-Y^2+t\in R[X,Y].
\end{equation}
The point $Q^\prime:=[\alpha :0:1]$ corresponds to a prime ideal
\[
q^\prime :=R[x,y](x-a)+R[x,y]y+R[x,y]t,\; a\in R^\times.
\]
Using equation (\ref{Z-chart}) we get
\[
t=y^2-s((x-a)+a)^3=y^2-s(x-a)((x-a)^2+3a(x-a)+3a^2)-sa^3 ,
\]
which shows that $t+sa^3\in q^\prime$. Since the element $t+sa^3$ generates $N$, we see that $q^\prime$ is generated by $x-a$ and $y$.

 %Inserted by TeXtelmExtel

We claim that $(x-a,y)$ is a regular sequence in $R[x,y]$, thus -- as performed earlier -- proving the regularity of $\cO_{\cC ,Q^\prime}$: indeed the isomorphism
\[
R[x,y]/(x-a)R[x,y]\cong R[Y]/(Y^2-(t+sa^3))
\]
shows that $x-a$ is a prime element of $R[x,y]$ and that $y$ is no zerodivisor on $R[x,y]/(x-a)R[x,y]$.

 %Inserted by TeXtelmExtel

Let $Q^{\prime\prime}:=[1:0:0]\in\cC_N$; so far we have shown that $\cC_N\setminus Q^{\prime\prime}\subset\Reg\cC$. The point $Q^{\prime\prime}$ itself is {\em not} regular: consider
\[
\cO_{\cC,Q^{\prime\prime}}=R[y,z]_{q^{\prime\prime}},\; q^{\prime\prime}:=R[y,z]y+R[y,z]z+R[y,z]t;
\]
if $\cO_{\cC,Q^{\prime\prime}}$ were regular, then according to the Theorems \ref{dimension equality} and \ref{maximal regular sequences} two among the three elements $y,z,t$ would be prime. However $t$ is not prime, because the local ring
\[
\cO_{\cC,Q^{\prime\prime}}/t\cO_{\cC,Q^{\prime\prime}}=\cO_{\cC_N,Q^{\prime\prime}}
\]
is no domain. In addition the factor ring
\[
\cO_{\cC,Q^{\prime\prime}}/z\cO_{\cC,Q^{\prime\prime}}=(R[y,z]/zR[y,z])_{q^{\prime\prime}}
\]
is no domain too.
\end{example}%
% ===================================
\section{Appendix: global dimension}
\label{s:bounds for gldim}
% ===================================
In general the finiteness of the weak dimension of a (local) ring $O$ has no influence on its global dimension. However Jensen \cite{Jen} and Osofsky \cite{Os2} give an upper bound for the global dimension in terms of the weak dimension provided that $O$ satisfies some cardinality condition. We use this generalization to obtain an upper bound for the global dimension $\gldim O$ of certain regular local rings $O\in\bL (S)$.

 %Inserted by TeXtelmExtel

Let $\kappa$ be an infinite cardinal number. A commutative ring $O$ is said to be {\bf $\kappa$-noetherian}  if every ideal $I\lhd O$ can be generated by a set of cardinality at most $\kappa$ and at least one ideal of $O$ cannot be generated by a set of cardinality less than $\kappa$. For every $m\in\mathbb{N}$ denote by $\aleph _m$ the $m$-th successor cardinal of the cardinality $\aleph_0$ of natural numbers. 
\begin{theorem}
\label{Global dimension}
Let $(S,N)$ be a non-trivial $\aleph _m$-noetherian valuation domain. Every essentially finitely generated regular $S$-algebra $O\in\bL (S)$ then satisfies $\gldim O\leq\wdim O+m+1\leq\dim\overline{O}+m+1$.
\end{theorem}
\proof
The second inequality follows from the first and Theorem \ref{upper bounds}.

 %Inserted by TeXtelmExtel

As for the first inequality we know from \cite{Jen}, Thm.~2 respectively \cite{Os2}, Cor.~1.4 that an $\aleph _m$-noetherian ring $O$ satisfies
\[
\gldim O\leq\wdim O+m+1,
\]
hence it remains to show that the local ring $O$ is $\aleph _l$-noetherian for some natural number $l\leq m$. 

 %Inserted by TeXtelmExtel

Since a localization of a $\kappa$-noetherian ring is $\kappa ^\prime$-noetherian for some $\kappa ^\prime\leq\kappa$, it suffices to prove that an algebra $A$ finitely generated over a $\aleph _m$-noetherian ring $S$ is $\aleph _{m^\prime}$-noetherian for some $m^\prime\leq m$. 

 %Inserted by TeXtelmExtel

As an $S$-module the algebra $A$ can be generated by countably many elements, for example by the monomials in a set $a_1,\ldots ,a_s$ of $S$-algebra generators of $A$. By \cite{Os2}, Lemma 1.1 this implies that every $S$-submodule of $A$ can be generated by $\aleph _m$ elements. Thus every ideal $I\lhd A$ as an $S$-module and hence as an $A$-module too can be generated by $\aleph _m$ elements, which completes the proof.
\eproof
\begin{footnotesize}

\end{footnotesize}
\end{document}